\documentstyle[A4,11pt]{article}

\begin{document}

\newcommand{\mod}{{\rm mod}\ }
\newcommand{\bbZ}{{\makebox{\bf{Z}}}}
\newtheorem{lemma}{Lemma}[section]
\newtheorem{theorem}{Theorem}[section]
\newtheorem{corollary}{Corollary}[section]
\newcommand{\bbox}{\hspace*{6mm}\rule{2mm}{3mm}}
\newcommand{\I}{\Pi}
\newcommand{\ol}{\overline}
\newcommand{\QED}{{\hspace*{\fill}$\Box$}}

\renewcommand{\thefootnote}{\fnsymbol{footnote}}

\title{On Determinants of Random\\
Symmetric Matrices over $\bbZ_m$} %

\author{Richard P. Brent\\
Computer Sciences Laboratory\\
Australian National University
\and
Brendan D.~McKay\\
Computer Science Department\\
Australian National University}

\date{Report TR-CS-88-03\\[10pt]
February 1988%
\footnotetext{\hspace*{-16pt}TR-CS-88-03 retyped by
Frances Page at Oxford University Computing
Laboratory, October 1998.\\
Copyright \copyright\ 1988--2010 the authors.
\hspace*{\fill} rpb101tr typeset using \LaTeX.}
}

\maketitle

\begin{abstract}
We determine the probability that a random $n$ by $n$ symmetric matrix
over $\{1, 2, \ldots, m\}$ has determinant divisible by $m$.
\end{abstract}

\thispagestyle{empty}

\section{Introduction}

Let $m$ be an integer.
The {\em m-rank\/} of an integer matrix $A$
is the greatest integer $k$ such that $A$ has a $k$ by $k$ submatrix
(not necessarily contiguous) whose determinant is nonzero $\mod m$,
or 0 if there is no such matrix.  If $m$ is a prime, the {\em m}-rank
is equivalent to the usual rank over the field $GF(m)$.  In this paper 
we assume that the elements $a_{ij}$ of $A$ are chosen at random,
independently and uniformly, from $\bbZ_m = \{1, 2, \ldots, m\}$,
subject to the condition that $a_{ij} = a_{ji}$, i.e., that $A$ is 
symmetric.  For corresponding results without the symmetry constraint,
see [1].

Let $P(n,m)$ denote the probability that a random $n$ by $n$ symmetric
matrix $A$ over $\bbZ_m$ has {\em m}-rank $n$, and define
$Q(n,m) = 1 - P(n,m)$.  Thus, $Q(n,m)$ is the probability that
\makebox{$\det(A) = 0 \;(\mod m)$.}
As in Lemma 1.1 of [1], we have

\begin{lemma}
Suppose $m = p^{\mu_1}_1 p^{\mu_2}_2 \cdots p^{\mu_k}_{k}$, where
$p_1, p_2, \ldots, p_k$ are distinctive primes.  Then

$$ Q(n,m) = \prod^k_{i=1}\; Q(n, p^{\mu_i}_i). $$ %
\end{lemma}

In view of Lemma 1.1, we restrict our attention to the case that $m$
is a prime power, say $m=p^\mu$.  It is useful to define $q=1/p$.

As in [1], our principal tool is Gaussian elimination, but in this case
we have to use forms of Gaussian elimination which preserve symmetry.
This is discussed in Section 2.  Then, in Section 3, we use symmetric
Gaussian elimination to show that $P(n,p^\mu)$ satisfies a five-term
recurrence relation.  In Section 4 we show that the five-term recurrence
can be reduced to a three-term recurrence.  Finally, in Section 5 we
show that the three-term recurrence can be solved explicitly.

The solution depends on the parity of $n$ and $\mu$, and is well-known
for $\mu = 1$, but appears to be new for $\mu > 1$.  To conclude
Section 5, we deduce some inequalities from the explicit solution.

An interesting problem is to determine the probability that a random
$n$ by $n$ symmetric matrix $A$ over $\bbZ_m$ has given {\em m}-rank
$r$, where $r<n$.  The case $\mu=1$ has been solved by Carlitz [2]
(provided $p \neq 2$), and the unsymmetric case has been considered in
[1], but the general symmetric case remains open.

Another open problem is to determine the probability that $\det(A)\; \mod m$
takes a given (nonzero) value $d$.  Small examples show that this 
probability depends on $d$.  For example, if $\mu=1, k=\lceil n/2
\rceil, (d|p)$ is the Legendre symbol,

$$
s = \left\{ \begin{array}{ll}
             0,      & {\rm\ if\ } p=2 {\rm\ or\ } n {\rm\ is\ odd},\\ 
             (d|p)(-1)^{k(p-1)/2}, & {\rm\ otherwise},\\
            \end{array} \right. $$

\medskip\noindent
and $\I$ is defined as in Section 5, then the probability
is

$$\left(\frac{q}{1-q}\right) \; \frac{\I_{2k}(q)}{\I_k(q^2)}\,(1+sq^k).$$

\section{Symmetric Gaussian Elimination}

Suppose that $A, n, m = p^\mu$ and $P(n,m)$ are as in Section 1.
Lemmas 2.1 and 2.2 describe symmetric versions of Gaussian elimination.

\begin{lemma}
Suppose that $n>1$ and $a_{11} \neq 0\; (\mod p)$.  Then there is a matrix
U such that

$$
U\!\!AU^T = \left(\begin{array}{c|ccc}
             a_{11} & 0 & \cdots & 0 \\ \hline
               0    &   &        &   \\
             \vdots &   &  A'    &    \\
               0    &   &        &
\end{array}\right)\; (\mod p^\mu)
\eqno{(2.1)}
$$

\noindent
where $A'$ is a random $n-1$ by $n-1$ symmetric matrix.
\end{lemma}

\noindent
{\bf Proof.}$\;$  Define

$$
      U = \left(\begin{array}{ccccc}
             1      &   &    &  & \\ 
              \lambda_2 & 1  &    & 0   &   \\
             \lambda_3 &   &  \ddots    &   & \\
               \vdots  &   &   0     & \ddots &  \\
             \lambda_n &    &       &       & 1 \\
\end{array}\right) 
$$

\medskip\noindent
where

$$a_{11}\lambda_j = -a_{1j} \; (\mod m) \;\;\; {\rm\ for\ }
j = 2,3,\ldots,n.
\eqno{(2.2)}$$

\medskip\noindent
Observe that $\lambda_2,\ldots,\lambda_n$ exist (since $a_{11} \neq
0 \;(\mod p)$) and (2.1) clearly holds.  Also, $A'$ depends linearly
on the random symmetric matrix
$$
            \left(\begin{array}{ccc}
             a_{22}& \cdots  &  a_{2n}    \\
               \vdots  & \ddots  &  \vdots   \\
             a_{n2} & \cdots    & a_{nn}
            \end{array} \right),
$$

\medskip\noindent
and thus is random over $\bbZ_m. $\QED

\begin{lemma}

Suppose that $n>2, a_{11} = 0\; (\mod p)$, and $a_{12} \neq 0\; (\mod p)$.
Then there is a matrix V such that

$$
V\!\!AV^T=   \left(\begin{array}{ccc|c}
             a_{11}\hspace*{-3.5mm}&  &\hspace*{-3.5mm}  a_{12}&    \\[-1.5ex]
                   &  &        & 0\\[-1.5ex]
             a_{21}\hspace*{-3.5mm}&  &\hspace*{-3.5mm} a_{22} &    \\\hline
                   &0&       & A''
            \end{array} \right)\;(\mod p^\mu)
\eqno{(2.3)}$$

\medskip\noindent
where $A''$ is a random $n-2$ by $n-2$ symmetric matrix.
\end{lemma}

\noindent
{\bf Proof.}$\;$ Define

$$
      V = \left(\begin{array}{ccccc}
             1      &   &    &     & \\ 
             0      & 1 &    & 0   &   \\
             \lambda_3 & \mu_3  & 1  &   & \\
               \vdots  & \vdots   &     & \ddots &  \\
             \lambda_n & \mu_n & 0  & &   1 \\
\end{array}\right),
$$

\medskip\noindent
where

$$
\pmatrix{ a_{11}&a_{12}\cr a_{21}&a_{22}\cr} 
\pmatrix{ \lambda_j \cr \mu_j}
= -\pmatrix{a_{1j} \cr a_{2j}} \;(\mod m) $$

\medskip\noindent
for $j=3,\ldots, n$.  Observe that $\mu_j$ and $\lambda_j$ 
exist $(j=3,\ldots,n)$ since

$$
\det \pmatrix{a_{11} & a_{12} \cr a_{21}& a_{22} \cr}
\neq 0 \;(\mod p).
$$

\medskip\noindent
It is easy to see that (2.3) holds.  Also, $A''$ depends linearly
on the random symmetric matrix

$$ \pmatrix{ 
a_{33} & \cdots & a_{3n} \cr
\vdots & \ddots & \vdots \cr
a_{n3} & \cdots & a_{nn} \cr } $$

\medskip\noindent
and thus is random over $\bbZ_m. $\QED

\section{A five-term recurrence for $P(n, p^\mu)$}

It is convenient to define

$$ P(0,p^\mu) = 1 \;\; {\rm\ for\ } \mu > 0 $$
\vspace*{-3mm}
\noindent
and \hfill (3.1)
\vspace*{-3mm}
$$ P(n,p^\mu) = 0 \;\; {\rm\ for\ } \mu \leq 0. $$

\begin{theorem}
If $n>0,\mu>0$, and boundary conditions are given by (3.1), then

$$P(n,p^\mu) = (1-q)P(n-1,p^\mu)+q(1-q^{n-1})P(n-2,p^\mu)$$
$$\hspace*{19mm}+\: q^n(1-q)P(n-1,p^{\mu-1})+q^{n+1}P(n,p^{\mu-2}).
\eqno{(3.2)}$$
\end{theorem}

\pagebreak[3]
{\samepage
\noindent
{\bf Proof.} $\;$ Let $A$ be a random symmetric $n \times n$ matrix
over $\bbZ_m, m=p^\mu$.  The four terms on the right side of (3.2)
arise from four mutually exclusive cases:

\begin{enumerate}
\item $a_{11} \neq 0 \; (\mod p)$.
\item $a_{11} = 0 \; (\mod p)$ and some $a_{1j} \neq 0 \; (\mod p)$.
\item $a_{1j} = 0 \; (\mod p)$ for $j=1,\ldots,n$ and
      $a_{11} \neq 0 \; (\mod p^2)$.
\item $a_{1j} = 0 \; (\mod p)$ for $j=1,\ldots,n$ and
      $a_{11} = 0 \; (\mod p^2)$.
\end{enumerate}
}

In case 1, which occurs with probability $1-q$, we apply Lemma 2.1;
since $\det(A)=a_{11}\det(A')$ $(\mod p^\mu)$ we have $\det(A) \neq 
0\;(\mod p^\mu)$ iff $\det(A') \neq 0\;(\mod p^\mu)$, which occurs
with conditional probability $P(n-1,p^\mu)$.

In case 2, which occurs with probability $q(1-q^{n-1})$, we can assume
that $a_{12}\neq 0\;(\mod p)$ by making a suitable permutation of rows
and columns, if necessary.  We then apply Lemma 2.2, obtaining
$\det(A)=(a_{11}a_{22}-a^2_{12})\det(A'')\;(\mod p^\mu)$, so
$\det(A)\neq0 \;(\mod p^\mu)$ iff %
$\det(A'')\neq 0$\linebreak[3] $(\mod p^\mu)$,
which occurs with conditional probability $P(n-2,p^\mu)$.

In case 3, which occurs with probability $q^n(1-q)$, we can reduce $A$
to the form (2.1) because (2.2) is solvable.  Thus 
$\det(A)=a_{11}\det(A')\;(\mod p^\mu)$ and
$\det(A) \neq 0 \;(\mod p^\mu)$ iff
$\det(A')\neq 0 \; (\mod p^{\mu-1})$,
which occurs with conditional probability $P(n-1,p^{\mu-1})$.

Finally, in case 4, which occurs with probability $q^{n+1}$, we can divide
the first row and column of $A$ by $p$, add random multiples of
$p^{\mu-1}$ to elements in the first row (and to corresponding elements
in the first column), add a random multiple of $p^{\mu-2}$ to the (1,1)
element, and obtain a new random symmetric matrix $\stackrel{\:\_}{A}$ such that

$$ \det(A)=p^2 \det(\stackrel{\:\_}{A}) \; (\mod p^\mu), $$

\medskip\noindent
so $\det(A)\neq0\; (\mod p^\mu)$ with conditional probability
$P(n,p^{\mu-2})$. \QED

\section{A three-term recurrence}

The five-term recurrence (3.2) with boundary conditions (3.1) can be used
to calculate $P(n,p^\mu)$ in $O(n\mu)$ arithmetic operations.
However, to obtain inequalities and asymptotic results it is useful 
to have an explicit solution.  To obtain such a solution, we first
reduce (3.2) to a three-term recurrence for $P(n,p^\mu)\;(n$ odd).

\begin{theorem}
If $\mu > 0$ and boundary conditions are given by (3.1), then for
odd $n \geq3$,

$$P(n,p^\mu)=(1-q^n)P(n-2,p^\mu)+q^n P(n,p^{\mu-2})
\eqno{(4.1)}$$
\noindent
and for odd $n\geq 1$,
$$P(n-1,p^\mu)= \frac{P(n,p^\mu)-q^n P(n,p^{\mu-1})}
{1-q^n}. \eqno{(4.2)}$$
\end{theorem}

\noindent
{\em Remarks.}$\;$ Equation (4.1) is a three-term recurrence from which
$P(n,p^\mu)$ can be calculated for odd~$n$.  Equation (4.2) then
determines $P(n,p^\mu)$ for even $n$.  Equations (4.1) and (4.2)
do not hold for all $n \geq 3$; for example, (4.2) fails if $n$ is
even and $\mu$ is odd.  In the unsymmetric case [1], (4.2) holds for both
even and odd $n$.\\

\pagebreak[3]
{\samepage
\noindent
{\bf Proof} (of Theorem 4.1) $\;$ Let

$$ {\cal P}_n = {\cal P}_n(x) = \sum^\infty_{\mu=1}P(n,p^\mu)x^\mu
\eqno{(4.3)}$$

\medskip\noindent
be a generating function for $P(n,p^\mu)$.  From the boundary conditions
(3.1) we have

$$ {\cal P}_0 = x/(1-x).  \eqno{(4.4)}$$

\medskip\noindent
Theorem 3.1 gives

$$ (1-q^{n+1}x^2){\cal P}_n = (1-q)(1+q^n x){\cal P}_{n-1} +
q(1-q^{n-1}){\cal P}_{n-2} \eqno{(4.5)}$$

\medskip\noindent
for $n\geq 1$. Thus, for odd $n=2k+1\geq 1$ we have

$$ (1-q^{2k+2}x^2){\cal P}_{2k+1} = (1-q)(1+q^{2k+1} x){\cal P}_{2k} +
q(1-q^{2k}){\cal P}_{2k-1} \eqno{(4.6)}$$

\medskip\noindent
and for even $n= 2k+2 \geq 2$ we have

$$ (1-q^{2k+3}x^2){\cal P}_{2k+2} = (1-q)(1+q^{2k+2} x){\cal P}_{2k+1} +
q(1-q^{2k+1}){\cal P}_{2k}. \eqno{(4.7)}$$

\medskip\noindent
Both (4.6) and (4.7) hold for $k \geq 0$, and with the boundary condition
(4.4) they define ${\cal P}_n$ for all $n \geq 0$.
}

Assume for the moment that (4.1) and (4.2) are correct.  From (4.1)
with $n=2k+3$ we have

$$
{\cal P}_{2k+3} = \left( \frac{1-q^{2k+3}}
{1-q^{2k+3}x^2} \right) {\cal P}_{2k+1}
\eqno{(4.8)}$$

\medskip\noindent
and from (4.2) with $n=2k+1$ we have

$$
{\cal P}_{2k} = \left( \frac{1-q^{2k+1}x}
{1-q^{2k+1}} \right) {\cal P}_{2k+1}.
\eqno{(4.9)}$$

\medskip\noindent
Clearly (4.8) and (4.9) for $k \geq 0$ and (4.4) define ${\cal P}_n$
for all $n \geq 0$.  Thus, it is sufficient to show that (4.8) and
(4.9) for $k \geq 0$ imply (4.6) and (4.7) for $k \geq 0$.  

From (4.8) and (4.9) we have

$$\hspace*{-4mm}{\cal P}_1 = \left( \frac{1-q}
{1-q x}\right) {\cal P}_0 \eqno{(4.10)} $$

$${\cal P}_2 = \left( \frac{1-q^3x}
{1-q^3 x^2}\right) {\cal P}_1, \eqno{(4.11)} $$

\medskip\noindent
which satisfy (4.6) and (4.7) with $k=0$.  Thus, we may assume
$k \geq 1$.  Equation (4.8) with $k$ replaced by $k-1$ gives 

$${\cal P}_{2k-1} = \left( \frac{1-q^{2k+1}x^2}
{1-q^{2k+1}}\right) {\cal P}_{2k+1}. \eqno{(4.12)} $$

\medskip\noindent
Substituting (4.9) and (4.12) in the right side of (4.8) and
simplifying, we obtain $(1-q^{2k+2}x^2){\cal P}_{2k+1}$, so (4.6)
holds.  Similarly, some algebra shows that (4.8) and (4.9) imply
(4.7).

Thus, the same generating function ${\cal P}_n, n \geq 0$, is
defined by (4.1) and (4.2) as by (4.5).  \QED\\

\pagebreak[3]
{\samepage
We can now give an explicit formula for the generating function
${\cal P}_n$.

\begin{theorem}
If $k \geq 0$, then

$${\cal P}_{2k+1}(x) =
\left( \frac{x}{1-x} \right)
\left( \frac{1-qx^2}{1-qx} \right)
\prod^k_{j=0}
\left( \frac{1-q^{2j+1}}{1-q^{2j+1}x^2} \right)
\eqno{(4.13)}$$
\noindent
and
$${\cal P}_{2k}(x) =
\left( \frac{1-q^{2k+1}x}{1-q^{2k+1}} \right) {\cal P}_{2k+1}(x).
\eqno{(4.14)}$$
\end{theorem}
}

\noindent
{\bf Proof.} $\;$ Equation (4.13) follows by induction from (4.8),
using (4.4) and (4.10).  Equation (4.14) is just (4.9). \QED

\section{An explicit solution and some bounds}

From Theorem 4.2 we can obtain an explicit solution for $P(n,p^\mu)$,
and hence for $Q(n,p^\mu)$. Define

$$ \I_n(q) = \prod^n_{j=1}\; (1-q^j)$$
\noindent
and
$$ T_\beta(k,s) = \left\{ \begin{array}{ll}
                   1, & {\rm if\ } k = 0,\\
                   {\displaystyle 
                    \sum^s_{j=0}q^{\beta j} \frac{\I_{k+j-1}(q^2)}
                             {\I_j (q^2)\I_{k-1} (q^2)}},&
                         {\rm if\ } k \geq 1.
                                              \end{array} \right.
\eqno{(5.1)}$$

\medskip\noindent
Our explicit solution may be written in terms of $T_1$ and $T_3$:

\begin{theorem}
If $n \geq 1, \mu \geq 1, k = \lfloor n/2\rfloor$, and $s = \lfloor(\mu 
- 1)/2\rfloor$ then

$$P(n, p^\mu) = \frac{\I_{2k}(q)}{(1-q)\I_k(q^2)}\,
\left( (1-q^{2k+1})T_3(k,s)-q^\mu(1-q^n)T_1(k,s )\right).
\eqno{(5.2)}$$
\end{theorem}

\medskip\noindent
{\bf Proof.}$\;$ We may show by induction on $k$ from (4.1) that

$$ P(2k + 1, p^\mu) = \frac{\I_{2k+1}(q)}{(1-q)\I_k(q^2)}\,
\left(T_3(k,s) - q^\mu T_1(k,s) \right).
\eqno{(5.3)}$$

\medskip\noindent
Considering the cases $\mu = 2s + 1$ and $\mu = 2s + 2$ separately, 
it follows from (4.2) that

$$P(2k,p^\mu) = P(2k + 1, p^\mu) + q^{2k + \mu}
\frac{\I_{2k}(q)}{\I_k(q^2)}\,T_1(k,s).
\eqno{(5.4)}$$

\medskip\noindent
After some simplification we see that (5.2) holds both for $n = 2k$
and $n = 2k+1$. \QED\\

Equation (5.2) is inconvenient for numerical computation as $P(n,p^\mu)$
is close to 1 unless $p^\mu$ is small.  In order to deduce a convenient
expression for $Q(n,p^\mu)= 1-P(n,p^\mu)$, we use the following
identities for $T_1$ and $T_3$.

\pagebreak[3]
{\samepage
\begin{lemma}
If $k \geq 0$ and $s \geq 0$ then

$$T_1(k,s)= \frac{\I_k(q^2)}{\I_{2k}(q)}\,
\left(1-q^{s+1}\sum^{k-1}_{j=0}q^{2j}
\frac{\I_{2j}(q)\I_{j+s}(q^2)}{\I_j(q^2)^2\I_s(q^2)}
\right).\eqno{(5.5)}$$

\medskip\noindent
and

$$T_3(k,s)= \frac{\I_k(q^2)}{\I_{2k+1}(q)}\,
\left(1-q-q^{3s+3}\sum^{k-1}_{j=0}q^{2j}
\frac{\I_{2j+1}(q)\I_{j+s}(q^2)}{\I_j(q^2)^2\I_s(q^2)}
\right).\eqno{(5.6)}$$
\end{lemma}
}

\noindent
{\bf Proof.}$\;$ (5.5) and (5.6) may be proved by induction on $k$.
The proof is similar to the proof of Theorem 2.1 of [1], so details
are omitted. \QED\\

If we substitute (5.5) into (5.2) the factor $\I_k(q^2)/\I_{2k}(q)$
cancels.  This gives a convenient explicit solution for $Q(n,p^\mu)$.

\begin{theorem}
If $n\geq 1, \mu \geq 1, k = \lfloor n/2 \rfloor$, and $s = \lfloor(\mu
- 1)/2 \rfloor$ then

$$
Q(n,p^\mu) = \frac{q^\mu(1-q^n)-R}{1-q}
\eqno{(5.7)}$$
\noindent where
$$R = q^{s+1}\sum^{k-1}_{j=0}\left(q^\mu(1-q^n)-q^{2s+2}(1-q^{2j+1})\right)
q^{2j}
\frac{\I_{2j}(q)\I_{j+s}(q^2)}{\I_j(q^2)^2\I_s(q^2)}
\eqno{(5.8)}$$
\noindent and
$$
0 \leq R < q^{3\mu/2}.
\eqno{(5.9)}$$
\end{theorem}

\noindent{\bf Proof.}$\;$ Equations (5.7) and (5.8) follow from Theorem
5.1 and Lemma 5.1 after some simplification.  Since $\mu \leq 2s + 2$
and $n \geq 2k$, we have $q^\mu(1-q^n)>q^{2s+2}(1-q^{2j+1})$ for
$0 \leq j \leq k-1$, so $R \geq 0$ (with equality only when $k=0$).
It is clear from (5.8) that $R=O(q^{\mu+s+1})=O(q^{3\mu/2})$, and
computation shows that $R<q^{3\mu/2}$ (the worst case is $n$ large,
$\mu = 2s+1$ large and $q=1/2$). \QED\\

We now show that $P(n,p^\mu)$ has the expected monotonicity properties.

\begin{theorem}
For $n \geq 0$ and $\mu \geq 0$,

$$P(n+1, p^\mu) \leq P(n,p^\mu)\leq P(n,p^{\mu + 1}).
\eqno{(5.10)}$$
\end{theorem}
\noindent{\bf Proof.}$\;$
The inequality $P(n,p^\mu) \leq P(n,p^{\mu+1})$ follows by induction 
on $n + \mu$ from the recurrence (3.2), since the coefficients on the
right side of (3.2) are independent of $\mu$.

To prove $P(n+1,p^\mu) \leq P(n,p^{\mu})$ we consider several cases.
If $n$ is even the inequality follows from (5.4).  If $n=2k+1$ is
odd and $\mu = 2s + 2$ is even then, from Theorem 5.1,

$$P(2k+1,p^{2s+2}) - P(2k+2,p^{2s+2})=q^{2k+2s+3}(1+q)
\frac{\I_{2k+1}(q)}{\I_k (q^2)}\; T_1(k+1,s) \geq 0.$$

\noindent
Finally, if $n=2k+1$ and $\mu = 2s+1$ are odd then, from Theorem 5.1,
{\samepage
\begin{eqnarray*}
\lefteqn{P(2k+1,p^{2s+1}) - P(2k+2,p^{2s+1})} \nonumber\\
&&= q^{2k+2s+2} \frac{\I_{2k+1}(q)}{\I_k(q^2)}
\left(T_1(k+1, s)-q^s\frac{\I_{k+s}(q^2)}{\I_s (q^2)}
\right)\nonumber\\
&&\geq 0. \nonumber\end{eqnarray*} \QED
}

\begin{corollary}
$\lim_{n \rightarrow \infty}Q(n,p^\mu)$ exists and lies in the interval
$[q^\mu, q^\mu/(1-q)]$.  Moreover,

$$
\lim_{n \rightarrow \infty} P(n,p^\mu)=
\frac{\I_\infty (q)}{(1-q)\I_\infty (q^2)}\,
\sum^s _{j=0}\frac{q^{3j}-q^{\mu+j}}{\I_j(q^2)}$$
\noindent and
$$
\lim_{n \rightarrow \infty} Q(n,p^\mu)=
\left(q^\mu-\frac{\I_\infty (q)}{\I_\infty (q^2)}\,
\sum^\infty_{j=s+1}\frac{q^{\mu + j}-q^{3j}}{\I_j(q^2)}
\right)\Biggl/(1-q),$$ %
\noindent where $s=\lfloor(\mu - 1)/2 \rfloor$.
\end{corollary}

\noindent
{\bf Proof.}$\;$ The limit exists by monotonicity in $n$, and the
bounds follow from this monotonicity and Theorem 5.2.  The explicit
limits follow from Theorem 5.1 and Lemma 5.1 respectively.  \QED

\section*{References}

[1] R.~P.~Brent and B.~D.~McKay,
Determinants and ranks of random matrices over \bbZ$_m$,
{\em Discrete Math.}~{\bf 66} (1987) 35-49.\\
\noindent[2] L.~Carlitz,
Representations by quadratic forms in a finite field,
{\em Duke Math.~J.}~{\bf 21} (1954) 123-137.

\end{document}